\documentclass[twoside]{article}
\usepackage{amsmath,amsthm,amssymb}
\usepackage{url}
\usepackage[spanish,english]{babel}
\usepackage{fancyhdr}
\fancypagestyle{plain}{\fancyhf{}%
\fancyhead[L]{\footnotesize{Divulgaciones Matem\'aticas Vol. 13 No. 2%
(2005), pp. \pageref{begin-art}--\pageref{end-art}}}%
}

\fancypagestyle{headings}{\fancyhf{}%
\fancyhead[LE,RO]{\thepage}%
\fancyhead[RE]{Ebner Pineda \& Wilfredo Urbina}%
\fancyhead[LO]{Non Tangential Convergence for the Ornstein-Uhlenbeck Semigroup.}%
\fancyfoot[C]{\footnotesize{Divulgaciones Matem\'aticas Vol. 13 No. 2%
(2005), pp. \pageref{begin-art}--\pageref{end-art}}}}

\theoremstyle{plain}
\newtheorem{thm}{Theorem}[section]

\newtheorem{lem}[thm]{Lemma}

\newtheorem{cor}[thm]{Corollary}

\theoremstyle{definition}
\newtheorem{defn}{Definition}[section]

\newcommand{\dem}{\medskip \par \noindent \mbox{\bf Proof. }} \def\ep{\hfill{$\Box $}}

\begin{document}
\label{begin-art}
\pagestyle{headings}
\thispagestyle{plain}
\footnote{\hspace{-18.1pt} Received ***. Revised ***.
Accepted ***. }
\selectlanguage{english}
\begin{center}
{\LARGE\bfseries Non Tangential Convergence for the Ornstein-Uhlenbeck Semigroup.\par }

\vspace{5mm}
{\large\slshape Convergencia no tangencial para el semigrupo de Ornstein-Uhlenbeck}\\[3mm]
{\large Ebner Pineda (\url{epineda@uicm.ucla.edu.ve})}\\[1mm]
Departamento de Matem\'{a}tica,  Decanato de Ciencia y Tecnolog\'{\i}a, UCLA\\
 Apartado 400 Barquisimeto 3001 Venezuela\\
\vspace{5mm}

{\large Wilfredo Urbina R.(\url{wurbina@euler.ciens.ucv.ve})}\\[1mm]
 Departamento de Matem\'{a}ticas, Facultad de Ciencias, UCV. \\
 Apartado 47195, Los Chaguaramos, Caracas 1041-A Venezuela, and\\
 Department of Mathematics and Statistics, University of New Mexico, \\
  Albuquerque, NM, 87131, USA.
\end{center}

\begin{abstract}
In this paper we are going to get the non tangential convergence,
in an appropriated parabolic  ``gaussian cone",  of the
Ornstein-Uhlenbeck semigroup in providing two proofs of this fact.
One is a direct proof by using the truncated non tangential
maximal function associated. The second one is obtained by using a
general statement. This second proof also allows us to get a
similar result for the Poisson-Hermite semigroup.

\end{abstract}

\selectlanguage{spanish}
\begin{abstract}
En este art\'{\i}culo vamos a obtener la convergencia no
tangencial, en un ``cono gaussiano " parab\'olico  apropiado, del
semigrupo de Ornstein-Uhlenbeck dando  dos pruebas diferentes de
ello. La primera es una prueba directa usando la funci\'on maximal
no tangencial truncada asociada. La segunda prueba se obtiene
usando principios generales. Esta \'ultima prueba nos permite
obtener un resultado an\'alogo para el semigrupo de
Poisson-Hermite.

{\bf Palabras y frases claves}: Convergencia no tangencial,
semigrupo de Ornstein-Uhlenbeck, semigrupo de Poisson-Hermite,
desarrollos de Hermite.
\end{abstract}

\selectlanguage{english}
\section{Introduction}

Let us consider the Gaussian measure $\gamma_d(x)=\frac{e^{-\left|x\right|^2}%
}{\pi^{d/2}}$ with $x\in\mathbb{R}^d$ and the Ornstein-Uhlenbeck differential
operator
\begin{equation}
L=\frac12\triangle_x-\left\langle x,\nabla _x\right\rangle.
\end{equation}

Let $\beta=(\beta _1,...,\beta_d)\in\mathbb{N}^d$ be a multi-index, let $\beta
!=\prod_{i=1}^d\beta _i!,$ $\left| \beta \right| =\sum_{i=1}^d\beta _i,$ $%
\partial _i=\frac \partial {\partial x_i},$ for each $1\leq i\leq d$ and $%
\partial ^\beta =\partial _1^{\beta _1}...\partial _d^{\beta _d}.$

Let us consider the normalized Hermite polynomial of order $\beta$, in $d$ variables
\begin{equation}
h_\beta (x)=\tfrac 1{\left( 2^{\left| \beta \right| }\beta !\right)
^{1/2}}\prod_{i=1}^d(-1)^{\beta _i}e^{x_i^2}\frac{\partial ^{\beta _i}}{%
\partial x_i^{\beta _i}}(e^{-x_i^2}),
\end{equation}
then, since the one dimensional Hermite polynomials satisfies the Hermite equation, see \cite{sz},  then the the normalized Hermite polynomial $h_\beta$ is an eigenfunction of $L$, with eigenvalue $-|\beta|$,
\begin{equation}\label{eigen}
L h_{\beta}(x)=-\left|\beta \right|h_\beta(x).
\end{equation}
Given a function $f$ $\in L^1(\gamma _d)$ its $\beta$-Fourier-Hermite
coefficient is defined by
\[
\hat{f}(\beta) =<f, h_\beta>_{\gamma_d} =\int_{\mathbb{R}^d}f(x)h_\beta (x)\gamma _d(dx).
\]
Let $C_n$ be the closed subspace of $L^2(\gamma_d)$ generated by
the linear combinations of $\left\{ h_\beta \ :\left| \beta
\right| =n\right\}$. By the orthogonality of the Hermite
polynomials with respect to $\gamma_d$ it is easy to see that
$\{C_n\}$ is an orthogonal decomposition of $L^2(\gamma_d)$,
$$ L^2(\gamma_d) = \bigoplus_{n=0}^{\infty} C_n$$
which is called the Wiener chaos.

Let $J_n$ be the orthogonal projection  of $L^2(\gamma_d)$ onto
$C_n$. If $f$ is a polynomial,
\[
J_n f=\sum_{\left|\beta\right|=n}\hat{f}(\beta) h_\beta.
\]
The Ornstein-Uhlenbeck semigroup $\left\{ T_t\right\}
_{t\geq 0}$ is given by
\begin{eqnarray}\label{01}
\nonumber T_t f(x)&=&\frac 1{\left( 1-e^{-2t}\right) ^{d/2}}\int_{\mathbb{R}^d}e^{-\frac{%
e^{-2t}(\left| x\right| ^2+\left| y\right| ^2)-2e^{-t}\left\langle
x,y\right\rangle }{1-e^{-2t}}}f(y)\gamma _d(dy)\\
& = & \frac{1}{\pi^{d/2}(1-e^{-2t})^{d/2}}\int_{\Bbb R^d}
e^{- \frac{|y-e^{-t}x|^2}{1-e^{-2t}}} f(y) dy.
\end{eqnarray}
$\left\{ T_t\right\}_{t\geq 0}$ is a strongly continuous Markov semigroup of contractions on $L^p(\gamma_d)$, with infinitesimal
generator $L$. Also, by a change of variable we can write,
\begin{equation}\label{t1}
T_t f(x)=\int_{\mathbb{R}^d} f(\sqrt{1-e^{-2t}}u + e^{-t}x)\gamma _d(du).
\end{equation}
\begin{defn} The maximal function for the Ornstein-Uhlenbeck semigroup is
defined as
\begin{eqnarray}
 T^*f(x) &=& \sup_{t > 0} |T_t f(x)| \nonumber \\
&=& \sup_{0<r<1}\frac{1}{\pi^{d/2}(1-r^2)^{d/2}}
|\int_{{\Bbb R}^d} e^{- \frac{|y-rx|^2}{1-r^2}} f(y) dy|.
\end{eqnarray}
\end{defn}

In \cite{gur} C. Guti\'errez and W. Urbina  obtained the following inequality for the
maximal function $ T^*f$,
\begin{equation} \label{des}
 T^* f(x) \leq C_d M_{\gamma_d} f(x) + (2 \vee |x|)^d e^{|x|^2} ||f||_{1,\gamma_d},
\end{equation}
 where $M_{\gamma_{d}}f$ is the Hardy-Littlewood maximal function of $f$ with respect to the gaussian measure $\gamma_{d}$,
 \begin{equation}\label{ineqmax}
M_{\gamma_d} f(x) = \sup_{r>0} \frac{1}{\gamma_d (B(x,r))} \int_{B(x,r)} |f(y)| \gamma_d(dy).
\end{equation}

Unfortunately, this inequality only allows to get the weak (1,1) continuity of $
T^*f$ in the  one dimensional case, $d=1$, but allows to get  a
pointwise convergence result. Several  results of this paper, see
Lemma 1.1 and Theorem 1.2, use techniques contained in that
paper.

If $f \in L^{1}(\gamma_{d}), u(x,t)=T_{t}f(x)$ is a solution of the
initial value problem
$$\left\{\begin{array}{cccl}
\dfrac{\partial u}{\partial t}(x,t)=L u(x,t)&&&\\
u(x,0)=f(x)&&&
\end{array}\right.$$
where $u(x,0)=f(x)$ means that
 $$\lim\limits_{t\to 0^{+}}u(x,t)=f(x), \, \, \mbox{a.e.} \, x$$
  We want to prove that this convergence is also non-tangential in the following sense.
  Let
\begin{equation}
\Gamma^{p}_{\gamma}(x)=\left\{(y,t)\in\mathbb{R}_{+}^{d+1}:\mid
y-x\mid<t^{\frac{1}{2}}\wedge\frac{1}{\mid x\mid}\wedge 1
\right\}
\end{equation}
be a parabolic  ``gaussian cone". We want  to prove that
\begin{equation*}
\lim_{(y,t)\rightarrow x,
\\(y,t)\in\Gamma^p_{\gamma}(x)}T_{t}f(y)=f(x), \, \, \mbox{a.e.}
\, x
\end{equation*}

Using the Bochner subordination formula (see \cite{se70}), 
\begin{equation*} 
 e^{-\lambda} = \frac{1}{\sqrt \pi} \int_0^{\infty} \frac{e^{-u}}{\sqrt u} e^{-\lambda^2/4u} du,
\end{equation*}
we define the Poisson-Hermite semigroup $\left\{ P_t\right\} _{t\geq 0}$ as
\begin{equation}\label{02}
P_t f(x)=\frac 1{\sqrt{\pi }}\int_0^{\infty} \frac{e^{-u}}{\sqrt{u}}T_{t^2/4u}f(x)du.
\end{equation}
$\left\{ P_t\right\}_{t\geq 0}$ is also a strongly continuous
semigroup on $L^p(\gamma_d)$, with infinitesimal generator
$-(-L)^{1/2}$. From (\ref{01}) we obtain, after the change of
variable $r=e^{-t^2/4u}$,
\begin{equation}\label{03}
P_t f(x)=\frac 1{2\pi ^{(d+1)/2}}\int_{\Bbb{R}^d}\int_0^1t\frac{\exp \left(
t^2/4\log r\right) }{(-\log r)^{3/2}}\frac{\exp \left( \frac{-\left|
y-rx\right| ^2}{1-r^2}\right) }{(1-r^2)^{d/2}}\frac{dr}rf(y)dy.
\end{equation}
\begin{defn} The maximal function for the Poisson-Hermite semigroup is
defined as
\begin{equation}
P^*f(x)  = \sup_{t > 0} |P_t f(x)|.
\end{equation}
\end{defn}
If $f \in L^{1}(\gamma_{d}), u(x,t)=P_{t}f(x)$ is solution of the
initial value problem
$$\left\{\begin{array}{cccl}
\dfrac{\partial^2 u}{\partial t^2}(x,t)=-L u(x,t)&&&\\
u(x,0)=f(x)&&&
\end{array}\right.$$
where $u(x,0)=f(x)$ means that
 $$\lim\limits_{t\to 0^{+}}u(x,t)=f(x),  \, \, \mbox{a.e.} \, x$$
 We want to prove that this convergence,  for the Poisson-Hermite
semigroup, is also non-tangential in the following sense.
  Let
\begin{equation}
 \Gamma_{\gamma}(x) = \left\{(y,t) \in \Bbb R^{d+1}_+ : \;\; |y-x| <
t \wedge \frac{1}{|x|} \wedge  1 \right\},
\end{equation}
 be a   ``gaussian cone".  Also we want to prove that
\begin{equation*}
\lim_{(y,t)\rightarrow x, \\(y,t)\in\Gamma_{\gamma}(x)}P_{t}f(y)=f(x),  \, \, \mbox{a.e.} \, x
\end{equation*}

In order to study the non-tangential convergence for the Ornstein-Uhlenbeck semigroup  we are going to consider the following maximal function, that was defined by L. Forzani and E. Fabes \cite{forfa}.

\begin{defn} The non tangential maximal function associated to the
Ornstein-Uhlenbeck semigroup is defined as
\begin{equation}
\mathcal{T}^{\ast}_{\gamma}f(x)=\sup_{(y,t)\in\Gamma^{p}_{\gamma}(x)}\mid
T_{t}f(y)\mid.
\end{equation}
\end{defn}

Using an inequality for a generalized maximal function, obtained by L. Forzani in \cite{for1} (for more details see \cite{uw98} pag 65--73 and 88--92),  it can be proved that $\mathcal{T}^{\ast}_{\gamma}f$ is weak $(1,1)$ and strong $(p,p)$ for $1<p< \infty$, with
respect to the Gaussian measure.

Actually for the non-tangential convergence for the Ornstein-Uhlenbeck semigroup it is enough to consider a ``truncated" maximal function. Let
\begin{equation}
 \Gamma^{p}(x)=\left\{(y,t) \in \Bbb R^{d+1}_+ : \;\; |y-x| <
t^{\frac{1}{2}},0<t< \frac{1}{|x|^{2}} \wedge  \frac{1}{4} \right\},
\end{equation}
be a truncated parabolic ``gaussian cone".
\begin{defn}The truncated non-tangencial maximal function
associated to the Ornstein-Uhlenbeck semigroup is defined as
\begin{equation}
\mathcal{T}^{\ast}f(x)=\sup_{(y,t)\in\Gamma^{p}(x)}\mid
T_{t}f(y)\mid.
\end{equation}
\end{defn}

\section{Non-tangential convergence of the
Ornstein-Uhlenbeck semigroup: direct proof.}

As we have mentioned already, the main result of this paper is to
prove the the non tangential convergence,
in an appropriated parabolic  ``gaussian cone",  of the
Ornstein-Uhlenbeck semigroup. Let us see 
 a direct proof by using the truncated non tangential
maximal function associated.

In the next lemma we are going to get a inequality better  than (\ref{ineqmax}) for the truncated non tangential maximal function $\mathcal{T}^{\ast}f$, which implies, immediately, that $\mathcal{T}^{\ast}f$ is weak $(1,1)$ and strong $(p,p)$ for $1<p< \infty$, with respect to the gaussian measure.
 
 \begin{lem}
\begin{equation}
\mathcal{T}^{\ast}f(x)\leq C_{d} M_{\gamma_{d}}f(x),
\end{equation}
for all $x\in\mathbb{R}^{d}$
\end{lem}

\dem
Let us take $u(y,t)=T_{t}f(y)$
and without loss of generality let us assume $f\geq 0$.

Let  $a_{o}=0$ and $a_{j}=\sqrt{j},\ \ j\in \mathbb{N},$ then
$a_{j}<a_{j+1}\ \ \ \forall j\in\mathbb{N}$, and let us denote
$$A_{j}(y,t)=\{u\in\mathbb{R}^{d}:a_{j-1}(1-e^{-2t})^{\frac{1}{2}}
\leq\mid e^{-t}y-u\mid<a_{j}(1-e^{-2t})^{\frac{1}{2}}\},$$
the annulus with center $e^{-t}y$. Now consider for each  $j\in \mathbb{N}$ the ball with center $e^{-t}y$,
and radius $a_{j}(1-e^{-2t})^{\frac{1}{2}}$ and let us denote it by $B_{j}(y,t)=B(e^{-t}y,a_{j}(1-e^{-2t})^{\frac{1}{2}}),$
then
$$A_{j}(y,t)=B_{j}(y,t)\setminus B_{j-1}(y,t)$$
$$\begin{array}{ccl}
u(y,t)&=&\dfrac{1}{\pi^{\frac{d}{2}}(1-e^{-2t})^{\frac{d}{2}}}
\displaystyle\int_{\mathbb{R}^{d}}e^{\dfrac{-\mid
e^{-t}y-u\mid^{2}}{1-e^{-2t}}}f(u)du\\
&=&\dfrac{1}{\pi^{\frac{d}{2}}(1-e^{-2t})^{\frac{d}{2}}}
\displaystyle\sum\limits_{j=1}^{\infty}
\int_{A_{j}(y,t)}e^{\dfrac{-\mid
e^{-t}y-u\mid^{2}}{1-e^{-2t}}}f(u)du
\end{array}$$
Now if $(y,t)\in\Gamma ^{p}(x)$ and $\mid
e^{-t}y-u\mid<a_{j}(1-e^{-2t})^{\frac{1}{2}}$ then,

$$\begin{array}{ccl}
\mid e^{-t}x-u\mid&=&\mid e^{-t}x-e^{-t}y+e^{-t}y-u\mid\\\\
&\leq&\mid e^{-t}(x-y)\mid+\mid e^{-t}y-u\mid\\\\
&<&e^{-t}t^{\frac{1}{2}}+a_{j}(1-e^{-2t})^{\frac{1}{2}}\\\\
&<&t^{\frac{1}{2}}+a_{j}(1-e^{-2t})^{\frac{1}{2}}\\\\
&<&(1+a_{j})(1-e^{-2t})^{\frac{1}{2}},
\end{array}$$
since $t<1-e^{-2t}$ if  $t<0.8$

Considering
$C_{j}(x,t)=B(e^{-t}x,(1+a_{j})(1-e^{-2t})^{\frac{1}{2}})$, we
have

$$u(y,t)\leq \dfrac{1}{\pi^{\frac{d}{2}}(1-e^{-2t})^{\frac{d}{2}}}
\sum_{j=1}^{\infty}e^{-a_{j-1}^{2}}\int_{C_{j}(x,t)}f(u)du$$

Now,
$$\begin{array}{l}
\displaystyle\int_{C_{j}(x,t)}f(u)du
=\displaystyle\int_{C_{j}(x,t)}f(u)e^{\mid
u\mid^{2}}e^{-\mid u\mid^{2}}du\\\\
=\displaystyle\int_{C_{j}(x,t)}f(u)e^{\mid
u-e^{-t}x\mid^{2}+2e^{-t}x.(u-e^{-t}x)+\mid e^{-t}x\mid^{2}}
e^{-\mid u\mid^{2}}du\\\\
\leq e^{(1+a_{j})^{2}(1-e^{-2t})+2(1+a_{j})(1-e^{-2t})^{\frac{1}{2}}\mid
e^{-t}x\mid+\mid e^{-t}x\mid^{2}} \displaystyle\int_{C_{j}(x,t)}f(u)e^{-\mid
u\mid^{2}}du
\end{array}$$

\noindent but, $\mid e^{-t}x-u\mid<(1+a_{j})(1-e^{-2t})^{\frac{1}{2}}$ and
therefore,
$$\mid
x-u\mid =\mid x-e^{-t}x+e^{-t}x-u\mid
<(1-e^{-t})\mid
x\mid+(1+a_{j})(1-e^{-2t})^{\frac{1}{2}}.$$

Taking $$D_{j}(x,t)=B(x,(1-e^{-t})\mid x\mid +(1+a_{j})
(1-e^{-2t})^{\frac{1}{2}}),$$

\noindent we get
$$\begin{array}{l}
\displaystyle\int_{C_{j}(x,t)}f(u)
e^{-\mid u\mid^{2}}du
\leq\displaystyle\int_{D_{j}(x,t)}f(u) e^{-\mid
u\mid^{2}}du\\\\
\leq M_{\gamma_{d}}f(x)\displaystyle\int_{D_{j}(x,t)} e^{-\mid
u\mid^{2}}du
=M_{\gamma_{d}}f(x)\displaystyle\int_{D_{j}(x,t)} e^{-\mid
u-x\mid^{2}+2x(x-u)-\mid x\mid^{2}}du\\\\
\leq M_{\gamma_{d}}f(x)e^{-\mid x\mid^{2}}\displaystyle\int_{D_{j}(x,t)}
e^{-\mid u-x\mid^{2}+2|x||x-u|}du\\\\
\leq M_{\gamma_{d}}f(x)e^{-\mid x\mid^{2}+2\mid
x\mid((1-e^{-t})\mid x\mid+(1+a_{j})(1-e^{-2t})^{\frac{1}{2}})}
\displaystyle\int_{D_{j}(x,t)}e^{-\mid u-x\mid^{2}}du\\\\
= M_{\gamma_{d}}f(x)e^{-\mid x\mid^{2}+2\mid x\mid((1-e^{-t})\mid
x\mid+(1+a_{j})(1-e^{-2t})^{\frac{1}{2}})}
\displaystyle\int_{E_{j}(x,t)}e^{-\mid w\mid^{2}}dw
\end{array}$$

\noindent where $E_{j}(x,t)=B(0,(1-e^{-t})\mid x\mid +(1+a_{j})
(1-e^{-2t})^{\frac{1}{2}})$. 

Since $\gamma_{d}$ is a $d-$dimensional measure, and using
that $t<\dfrac{1}{\mid x\mid^{2}}\wedge\frac{1}{4}$, we get

\begin{eqnarray*}
\displaystyle\int_{C_{j}(x,t)}f(u)e^{-\mid
u\mid^{2}}du &\leq& C_{d} M_{\gamma_{d}}f(x)e^{-\mid x\mid^{2}+2\mid
x\mid((1-e^{-t})\mid x\mid+(1+a_{j})(1-e^{-2t})^{\frac{1}{2}})}\\
&&\quad \quad  \quad \quad \quad \times ((1-e^{-t})\mid x\mid+(1+a_{j})(1-e^{-2t})^{\frac{1}{2}})^{d}\\
&=&C_{d}M_{\gamma_{d}}f(x)e^{-\mid x\mid^{2}+2\mid
x\mid((1-e^{-t})\mid x\mid+(1+a_{j})(1-e^{-2t})^{\frac{1}{2}}}\\
&&  \quad \quad \times(1-e^{-t})^{\frac{d}{2}}((1-e^{-t})^{\frac{1}{2}}\mid x\mid+(1+a_{j})(1+e^{-t})^{\frac{1}{2}})^{d}\\
&\leq& C_{d}M_{\gamma_{d}}f(x)e^{-\mid
x\mid^{2}+2\frac{(1-e^{-t})}{t}+2(1+a_{j})\frac{(1-e^{-2t})^{\frac{1}{2}}}
{t^{\frac{1}{2}}}}\\
&& \quad \quad \times(1-e^{-t})^{\frac{d}{2}}\left(\dfrac{(1-e^{-t})^{\frac{1}{2}}}{t^{\frac{1}{2}}}
+(1+a_{j})(1+e^{-t})^{\frac{1}{2}}\right)^{d}.
\end{eqnarray*}
 Therefore
$$\begin{array}{l}
\displaystyle\int_{C_{j}(x,t)}f(u)du
\leq
e^{(1+a_{j})^{2}(1-e^{-2t})+2(1+a_{j})\frac{(1-e^{-2t})^{\frac{1}{2}}}
{t^{\frac{1}{2}}}+e^{-2t}\mid x\mid^{2}}\displaystyle\int_{C_{j}(x,t)}f(u)
e^{-\mid u\mid^{2}} du\\\\
\leq
e^{(1+a_{j})^{2}(1-e^{-2t})+2(1+a_{j})\frac{(1-e^{-2t})^{\frac{1}{2}}}
{t^{\frac{1}{2}}}+\mid x\mid^{2}}
C_{d}M_{\gamma_{d}}f(x)e^{-|x|^{2}+2\frac{(1-e^{-t})}{t}+2(1+a_{j})
\frac{(1-e^{-2t})^{\frac{1}{2}}}
{t^{\frac{1}{2}}}}\\\hspace{3cm}
\times (1-e^{-t})^{\frac{d}{2}}(\dfrac{(1-e^{-t})^{\frac{1}{2}}}{t^{\frac{1}{2}}}
+(1+a_{j})(1+e^{-t})^{\frac{1}{2}})^{d}\\\\
\leq
e^{(1+a_{j})^{2}(1-e^{-\frac{1}{2}})+4(1+a_{j})\frac{(1-e^{-2t})^{\frac{1}{2}}}
{t^{\frac{1}{2}}}+\frac{2(1-e^{-t})}{t}}(1-e^{-t})^{\frac{d}{2}}\\\hspace{3cm}
\times
\left(\dfrac{(1-e^{-2t})^{\frac{1}{2}}}{t^{\frac{1}{2}}}+(1+a_{j})
(1+e^{-t})^{\frac{1}{2}}\right)^{d}C_{d}M_{\gamma_{d}}f(x)\\
\leq e^{(1+a_{j})^{2}(1-e^{-\frac{1}{2}})+4(1+a_{j})\sqrt{2}+2}.
(1-e^{-t})^{\frac{d}{2}}(1+(1+a_{j})\sqrt{2})^{d}C_{d}M_{\gamma_{d}}f(x),
\end{array}$$
since $0<t<\frac{1}{4}$ and $\dfrac{1-e^{-t}}{t}<1,1+e^{-t}<2$, if
$t>0$.

Thus,
\begin{eqnarray*}
u(y,t)&\leq& \dfrac{1}{\pi^{\frac{d}{2}}(1-e^{-2t})^{\frac{d}{2}}}
\displaystyle\sum_{j=1}^{\infty}e^{-a_{j-1}^{2}}\displaystyle\int_{C_{j}(x,t)}f(u)du\\
&\leq& C_{d}M_{\gamma_{d}}f(x)
\dfrac{1}{\pi^{\frac{d}{2}}(1+e^{-t})^{\frac{d}{2}}
(1-e^{-t})^{\frac{d}{2}}}\\
&&  \times \displaystyle\sum_{j=1}^{\infty}e^{-a_{j-1}^{2}}
e^{(1+a_{j})^{2}(1-e^{-\frac{1}{2}})+4(1+a_{j})\sqrt{2}+2}
 (1-e^{-t})^{\frac{d}{2}}(1+(1+a_{j})\sqrt{2})^{d}\\
&\leq& C_{d}M_{\gamma_{d}}f(x)\dfrac{1}{\pi^{\frac{1}{2}}}
\displaystyle\sum_{j=1}^{\infty}e^{-a_{j-1}^{2}+(1+a_{j})^{2}
(1-e^{-\frac{1}{2}})+4(1+a_{j})\sqrt{2}+2}(1+(1+a_{j})\sqrt{2})^{d},
\end{eqnarray*}
since $1+e^{-t}\geq 1.$ Now it is easy to see that
\begin{eqnarray*}
&&-a_{j-1}^{2}+(1+a_{j})^{2}(1-e^{-\frac{1}{2}})+4(1+a_{j})\sqrt{2}+2\\
&&\quad \quad =4+4\sqrt{2}-e^{-\frac{1}{2}}-[-(2(1-e^{-\frac{1}{2}})+4\sqrt{2})+
e^{-\frac{1}{2}}\sqrt{j}]\sqrt{j},
\end{eqnarray*}

which is negative for $j$ sufficiently big, then
$$\sum_{j=1}^{\infty}e^{-a_{j-1}^{2}+(1+a_{j})^{2}
(1-e^{-\frac{1}{2}})+4(1+a_{j})\sqrt{2}+2}.(1+(1+a_{j}).\sqrt{2})^{d}
<\infty.$$
Thus $u(y,t)\leq C_{d}M_{\gamma_{d}}f(x)$ and since $(y,t)\in \Gamma^{p}(x)$ is arbitrary
$$\mathcal{T}^{\ast}f(x)=\sup_{(y,t)\in\Gamma^{p}(x)}u(y,t)\leq C_{d}
M_{\gamma_{d}}f(x). $$ \ep

Now we are ready to establish the convergence result for the  Ornstein-Uhlenbeck semigroup.

\begin{thm} The Ornstein-Uhlenbeck semigroup
$\{T_{t}f\}$ converges in $L^1(\gamma_{d})$ a.e if $t\to
0^{+}$, for any function $f\in L^1(\gamma_{d})$,
\begin{equation}
\lim\limits_{t\to 0^{+}}u(x,t)=f(x), \, \, \mbox{a.e.} \, x
\end{equation}
Moreover, if $u(y,t)=T_{t}f(y)$ then $u(y,t)$ tends to $f(x)$ non
tangentially ,i.e.
\begin{equation}
\lim_{(y,t)\rightarrow x,
\\(y,t)\in\Gamma^p_{\gamma}(x)}T_{t}f(y)=f(x), \, \, \mbox{a.e.}
\, x.
\end{equation}
\end{thm}

\dem
We have,
$$u(y,t)=\dfrac{1}{\pi^{\frac{d}{2}}(1-e^{-2t})^{\frac{d}{2}}}
\int_{\mathbb{R}^{d}}e^{\dfrac{-\mid
e^{-t}y-u\mid^{2}}{1-e^{-2t}}}f(u)du,$$
considering
$$\Omega f(x)=\lim_{\alpha\to 0^{+}}\left[\sup_{(y,s)\in\Gamma^{p}_{\gamma}
(x),0<s<\alpha}\mid u(y,s)-f(x)\mid\right],$$
and let us set
$f(x)=f(x)\chi_{(0,k)}+f(x)(I-\chi_{(0,k)})=f_{1}(x)+f_{2}(x),$ for
$k\in\mathbb{N}$ fix.

Let us prove that
$$\Omega f(x)\leq C_{d}M_{\gamma_{d}}f_{2}(x), \mbox{a.e},$$
for $\mid x\mid\leq k-1$.

Let us consider  $x$ a
Lebesgue's point for $f\in L^1(\gamma_{d})$, i.e. $x$ verifies
$$\lim_{r\to 0^{+}}\dfrac{1}{\gamma_{d}(B(x;r))}\int_{B(x;r)}\mid
f(u)-f(x)\mid\gamma_{d}(du)=0 \,\,\,$$

Then given  $\epsilon>0$ there exists
$0<\delta<1$ such that
$$\dfrac{1}{\gamma_{d}(B(x;r))}\int_{B(x;r)}\mid
f(u)-f(x)\mid\gamma_{d}(du)<\epsilon,$$
for $0<r<\delta.$ Let us define $g$ as
$g(u)=\left\{\begin{array}{lcl}
f(u)-f(x)&if&\mid u-x\mid\leq\delta\\
0&if&\mid u-x\mid>\delta \end{array}\right.$
Thus $g$ depends
on $x$ and  $M_{\gamma_{d}}g(x)<\epsilon$.

On the other hand, since $$u(y,t)-f(x)=u^1(y,t)-f_1(x)+u^{2}(y,t)-f_{2}(x)$$
where
$$u^{i}(y,t)=\dfrac{1}{\pi^{\frac{d}{2}}(1-e^{-2t})^{\frac{d}{2}}}
\int_{\mathbb{R}^{d}}e^{\dfrac{-\mid
e^{-t}y-u\mid^{2}}{1-e^{-2t}}}f_i(u)du\  \ i=1,2,$$
then we get,
$$\begin{array}{l}
u^{1}(y,t)-f_1(x)
=\dfrac{1}{\pi^{\frac{d}{2}}(1-e^{-2t})^{\frac{d}{2}}}
 \displaystyle\int_{\mathbb{R}^{d}}e^{\dfrac{-\mid
e^{-t}y-u\mid^{2}}{1-e^{-2t}}}(f_{1}(u)-f_{1}(x))du\\\\
=\dfrac{1}{\pi^{\frac{d}{2}}(1-e^{-2t})^{\frac{d}{2}}}
 \displaystyle\int_{\mid x-u\mid\leq\delta}e^{\dfrac{-\mid
e^{-t}y-u\mid^{2}}{1-e^{-2t}}}(f_{1}(u)-f_{1}(x))du  \\\\ +
\dfrac{1}{\pi^{\frac{d}{2}}(1-e^{-2t})^{\frac{d}{2}}}  \displaystyle\int_{\mid
x-u\mid>\delta}e^{\dfrac{-\mid
e^{-t}y-u\mid^{2}}{1-e^{-2t}}}(f_{1}(u)-f_{1}(x))du.
\end{array}$$
Now we have that if $\mid x\mid\leq k-1$ and $(y,t)\in\Gamma ^{p}_{\gamma}(x)$
with $t< \frac{1}{|x|^{2}} \wedge  \frac{1}{4}$, then
$(y,t)\in\Gamma ^{p}(x)$. Thus $\mid u-x\mid\leq\delta$ implies
$$\mid u\mid=\mid u-x+x\mid\leq \mid
u-x\mid+\mid x\mid<\delta+k-1<1+k-1=k$$
and then, $f_{1}(u)=f(u)\wedge f_{1}(x)=f(x).$ Therefore
$$\begin{array}{ccl}
&&  \dfrac{1}{\pi^{\frac{d}{2}}(1-e^{-2t})^{\frac{d}{2}}}
 \left| \displaystyle\int_{\mid x-u\mid\leq\delta}e^{\dfrac{-\mid
e^{-t}y-u\mid^{2}}{1-e^{-2t}}}(f_{1}(u)-f_{1}(x))du  \right| \\
&=&\dfrac{1}{\pi^{\frac{d}{2}}(1-e^{-2t})^{\frac{d}{2}}} \left|\displaystyle\int_{\mid
x-u\mid\leq\delta}e^{\dfrac{-\mid
e^{-t}y-u\mid^{2}}{1-e^{-2t}}}(f(u)-f(x))du\right|\\\\
&=& \dfrac{1}{\pi^{\frac{d}{2}}(1-e^{-2t})^{\frac{d}{2}}} \left|
\displaystyle\int_{\mathbb{R}^{d}}e^{\dfrac{-\mid
e^{-t}y-u\mid^{2}}{1-e^{-2t}}} g(u) du\right|\\\\
&\leq&\mathcal{T}^{\ast}g(x)
\leq C_{d}M_{\gamma_{d}}g(x)
\leq C_{d}\epsilon.
\end{array}$$
Now observe that if $(y,t)\in\Gamma ^{p}_{\gamma}(x)$ and
$t^{\frac{1}{2}}\leq\frac{\delta}{2}$ then, $\mid u-x\mid>\delta$
implies $\delta<\mid u-x\mid\leq\mid u-y\mid+\mid y-x\mid$ and
thus
$$\delta<\mid u-y\mid+\mid y-x\mid<\mid
u-y\mid+t^{\frac{1}{2}}\leq\mid u-y\mid+\frac{\delta}{2},$$
thus $\mid u-y\mid>\frac{\delta}{2}.$ Therefore,
\begin{eqnarray*}
&&
\dfrac{1}{\pi^{\frac{d}{2}}(1-e^{-2t})^{\frac{d}{2}}}
\left|\displaystyle\int_{\mid u-x\mid>\delta}e^{\dfrac{-\mid e^{-t}y-u\mid^{2}}{1-e^{-2t}}}(f_{1}(u)-f_{1}(x))du\right|\\
\end{eqnarray*}
\begin{eqnarray*}
&\leq&\dfrac{1}{\pi^{\frac{d}{2}}(1-e^{-2t})^{\frac{d}{2}}}\displaystyle\int_{\mid u-x\mid>\delta}
e^{\dfrac{-\mid e^{-t}y-u\mid^{2}}{1-e^{-2t}}} \mid f_{1}(u)\mid du\\
& &\quad \quad \quad \quad + \dfrac{1}{\pi^{\frac{d}{2}}(1-e^{-2t})^{\frac{d}{2}}}\mid f_{1}(x)\mid
\displaystyle \int_{\mid u-x\mid>\delta}e^{\dfrac{-\mid e^{-t}y-u\mid^{2}}{1-e^{-2t}}}du\\
&\leq&\dfrac{1}{\pi^{\frac{d}{2}}(1-e^{-2t})^{\frac{d}{2}}}
 \displaystyle\int_{\mid u-y\mid>\frac{\delta}{2}} e^{\dfrac{-\mid e^{-t}y-u\mid^{2}}{1-e^{-2t}}}
 \mid f_{1}(u)\mid du\\
& &\quad \quad \quad \quad + \dfrac{1}{\pi^{\frac{d}{2}}(1-e^{-2t})^{\frac{d}{2}}} \mid
f_{1}(x)\mid \int_{\mid u-x\mid>\delta}e^{\dfrac{-\mid e^{-t}y-u\mid^{2}}{1-e^{-2t}}}du.
\end{eqnarray*}

Now, we have
$$\begin{array}{ccl}
&& \dfrac{1}{\pi^{\frac{d}{2}}(1-e^{-2t})^{\frac{d}{2}}}
 \displaystyle\int_{\mid u-y\mid>\frac{\delta}{2}} e^{\dfrac{-\mid e^{-t}y-u\mid^{2}}{1-e^{-2t}}}
 \mid f_{1}(u)\mid du\\
&=&\dfrac{1}{\pi^{\frac{d}{2}}(1-e^{-2t})^{\frac{d}{2}}}
\displaystyle\int_{\mid u-y\mid>\frac{\delta}{2},\mid u\mid<k}e^{\dfrac{-\mid
e^{-t}y-u\mid^{2}}{1-e^{-2t}}}\mid f_{1}(u)\mid du\\\\
&=&\dfrac{1}{\pi^{\frac{d}{2}}(1-e^{-2t})^{\frac{d}{2}}}
\displaystyle\int_{\mid u-y\mid>\frac{\delta}{2},\mid u\mid<k}e^{\dfrac{-\mid
e^{-t}y-u\mid^{2}}{1-e^{-2t}}}\mid f(u)\mid du\\\\
&=&\dfrac{1}{\pi^{\frac{d}{2}}(1-e^{-2t})^{\frac{d}{2}}}
\displaystyle\int_{\mid u-y\mid>\frac{\delta}{2},\mid u\mid<k}e^{\dfrac{-\mid
e^{-t}y-u\mid^{2}}{1-e^{-2t}}}e^{\mid u\mid^{2}}\mid f(u)\mid
e^{-\mid u\mid^{2}}du\\\\
\end{array}$$
$$\begin{array}{ccl}
&\leq&\dfrac{1}{\pi^{\frac{d}{2}}(1-e^{-2t})^{\frac{d}{2}}}e^{k^{2}}
\displaystyle\int_{\mid u-y\mid>\frac{\delta}{2},\mid u\mid<k}e^{\dfrac{-\mid
e^{-t}y-u\mid^{2}}{1-e^{-2t}}}\mid f(u)\mid e^{-\mid u\mid^{2}}du.
\end{array}$$
Then for $0<t<\log\left(\dfrac{4k+2\delta}{4k+\delta}\right);\mid
u-y\mid>\frac{\delta}{2},\mid u\mid<k$ implies that

$$\begin{array}{ccl}
\mid e^{-t}y-u\mid&=&\mid e^{-t}y-e^{-t}u+e^{-t}u-u\mid=\mid
e^{-t}(y-u)-(u-e^{-t}u)\mid \\\\
&\geq&e^{-t}\mid y-u\mid -\mid u-e^{-t}u\mid=e^{-t}\mid y-u\mid
-(1-e^{-t})\mid u\mid\\\\
&\geq&e^{-t}\dfrac{\delta}{2}-k(1-e^{-t})=e^{-t}\left(\dfrac{\delta}{2}+k
\right)-k,
\end{array}$$
but $0<t<\log\left(\dfrac{4k+2\delta}{4k+\delta}\right)$ and
therefore $e^{-t}>\dfrac{4k+\delta}{4k+2\delta},$ then,
$$\begin{array}{ccl}
e^{-t}\left(\dfrac{\delta}{2}+k\right)-k&>&\dfrac{4k+\delta}{4k+2\delta}
\left(\dfrac{\delta+2k}{2}\right)-k\\\\
&=&\dfrac{4k+\delta}{4(2k+\delta)}(2k+\delta)-k\\\\
&=&\dfrac{4k+\delta}{4}-k
=\dfrac{4k+\delta-4k}{4}=\dfrac{\delta}{4}.
\end{array}$$

\noindent Therefore $\mid u-y\mid >\dfrac{\delta}{2},\ \mid u\mid< k$ implies $\mid e^{-t}y-u\mid>\dfrac{\delta}{4}$  and thus
$$\begin{array}{ccl}
&&\dfrac{1}{\pi^{\frac{d}{2}}(1-e^{-2t})^{\frac{d}{2}}}
 \displaystyle\int_{\mid u-y\mid>\frac{\delta}{2}} e^{\dfrac{-\mid e^{-t}y-u\mid^{2}}{1-e^{-2t}}}
 \mid f_{1}(u)\mid du\\
 &\leq&\dfrac{1}{\pi^{\frac{d}{2}}(1-e^{-2t})^{\frac{d}{2}}}e^{k^{2}}
\displaystyle\int_{\mid u-y\mid>\frac{\delta}{2},\mid u\mid<k}e^{-\frac{\delta^{2}}
{16(1-e^{-2t})}}\mid f(u)\mid e^{-\mid u\mid{}^{2}}du\\\\
&\leq&\dfrac{e^{-\frac{\delta^{2}}{16(1-e^{-2t})}+k^{2}}}{\pi^{\frac{d}{2}}
(1-e^{-2t})^{\frac{d}{2}}}\displaystyle\int_{\mathbb{R}^{d}}\mid f(u)\mid
e^{-\mid u\mid{}^{2}}du=\dfrac{e^{-\frac{\delta^{2}}{16(1-e^{-2t})}+k^{2}}}{\pi^{\frac{d}{2}}
(1-e^{-2t})^{\frac{d}{2}}}\|f\|_{1,\gamma_{d}}
\end{array}$$

On the other hand, taking the change of variable $s=u-e^{-t}y$, we have 
$$\begin{array}{ccl}
&& \dfrac{1}{\pi^{\frac{d}{2}}(1-e^{-2t})^{\frac{d}{2}}} \mid
f_{1}(x)\mid \int_{\mid u-x\mid>\delta}e^{\dfrac{-\mid e^{-t}y-u\mid^{2}}{1-e^{-2t}}}du\\\\
&=&\dfrac{\mid f_{1}(x)\mid}{\pi^{\frac{d}{2}}(1-e^{-2t})^{\frac{d}{2}}}
\displaystyle\int_{\mid x-s-e^{-t}y\mid>\delta} e^{\frac{-\mid s\mid^{2}}{1-e^{-2t}}}ds\\\\
&=&\dfrac{\mid f(x)\mid}{\pi^{\frac{d}{2}}(1-e^{-2t})^{\frac{d}{2}}}
\displaystyle\int_{\mid x-s- e^{-t}y\mid>\delta} e^{\frac{-\mid s\mid^{2}}{1-e^{-2t}}}ds,
\end{array}$$
since, $f_{1}(x)=f(x)$ as $\mid x\mid \leq k-1< k.$

Thus taking $0<t<\log\left(\dfrac{k-1-\delta/2}{k-1-3\delta/4}\right),$ $\mid x-s-e^{-t}y\mid>\delta$
implies
$$\begin{array}{ccl}
 \mid s\mid&=&
\mid s-x+e^{-t}y+x-e^{-t}y\mid
=\mid s-x+e^{-t}y-(e^{-t}y-x)\mid\\\\
&\geq&\mid s-x+e^{-t}y\mid-\mid e^{-t}y-x\mid.
\end{array}$$

But
$$\begin{array}{ccl}
\mid e^{-t}y-x\mid&=&\mid e^{-t}y-e^{-t}x+e^{-t}x-x\mid
\leq e^{-t}\mid y-x\mid +(1-e^{-t})\mid x\mid\\\\
&\leq&e^{-t}t^{\frac{1}{2}}+(1-e^{-t})(k-1).
\end{array}$$
Thus, since $t^{\frac{1}{2}}\leq \dfrac{\delta}{2}$,
$$\begin{array}{ccl}
\mid s-x+e^{-t}y\mid-\mid e^{-t}y-x\mid&>& \delta-e^{-t}
t^{\frac{1}{2}}-(1-e^{-t})(k-1)\\\\
&\geq&\delta-e^{-t}\dfrac{\delta}{2}-(k-1)(1-e^{-t})\\\\
&=&\delta-(k-1)+(k-1-\dfrac{\delta}{2})e^{-t},
\end{array}$$

\noindent and as $0<t<\log\left(\dfrac{k-1-\delta/2}{k-1-3\delta/4}\right),$
then $e^{-t}>\dfrac{k-1-3\delta/4}{k-1-\delta/2}$. Hence,
$$\begin{array}{ccl}
\mid s\mid&>&\delta-(k-1)+(k-1-\delta/2)e^{-t}\\\\
&>&\delta-(k-1)+(k-1-\delta/2)\dfrac{k-1-3\delta/4}{k-1-\delta/2}\\\\
&=&\delta-(k-1)+k-1-3\delta/4=\delta-3\delta/4=\dfrac{\delta}{4}.
\end{array}$$

Then
$\mid x-s-e^{-t}y\mid >\delta$ implies $ \mid s\mid>\dfrac{\delta}{4}$ if $
0<t<\log\left(\dfrac{k-1-\delta/2}{k-1-3\delta/4}\right)$.
Therefore, taking $w=\dfrac{s}
{\sqrt{1-e^{-2t}}},$
$$\begin{array}{ccl}
&& \dfrac{1}{\pi^{\frac{d}{2}}(1-e^{-2t})^{\frac{d}{2}}}
 \displaystyle\int_{\mid u-y\mid>\frac{\delta}{2}} e^{\dfrac{-\mid e^{-t}y-u\mid^{2}}{1-e^{-2t}}}
 \mid f_{1}(u)\mid du\\
&\leq&\dfrac{\mid f(x)\mid}{\pi^{\frac{d}{2}}(1-e^{-2t})^{\frac{d}{2}}}
\displaystyle\int_{\mid s\mid>\frac{\delta}{4}}e^{\dfrac{-\mid s\mid{}^{2}}{1-e^{-2t}}}ds\\
&=&\dfrac{\mid f(x)\mid}{\pi^{\frac{d}{2}}}\displaystyle\int_{\mid w\mid >\frac{\delta}
{4\sqrt{1-e^{-2t}}}}e^{-\mid w\mid{}^{2}} dw.
\end{array}$$

Now since, $\mid x\mid\leq k-1<k$, then$f_{2}(x)=0$. Hence
$$\mid u^{2}(y,t)-f_{2}(x)\mid =\mid u^{2}(y,t)\mid\leq \mathcal{T}^{\ast}f_{2}(x)\leq
C_{d}M_{\gamma_{d}}f_{2}(x)$$
for $(y,t)\in\Gamma ^{p}(x).$ Therefore,
$$\begin{array}{l}
\mid u(y,t)-f(x)\mid \leq \mid u^{1}(y,t)-f_{1}(x)\mid +\mid
u^{2}(y,t)-f_{2}(x)\mid\\\\
=\mid u^{1}(y,t)-f_{1}(x)\mid+\mid u^{2}(y,t)\mid \\\\
\leq C_{d}\epsilon+\dfrac{e^{\frac{-\delta}{16(1-e^{-2t})}+k^{2}}}
{(1-e^{-2t})^{\frac{d}{2}}}\|f\|_{1,\gamma_{d}}+\dfrac{\mid f(x)\mid}
{\pi^{\frac{d}{2}}} \displaystyle\int_{\mid w\mid >\frac{\delta}{4\sqrt{1-e^{-2t}}}}
e^{-\mid w\mid{}^{2}} dw \\\\ \quad \quad \quad \quad \quad \quad \quad \quad \quad \quad \quad \quad \quad \quad \quad \quad +C_{d}M_{\gamma_{d}}f_{2}(x),
\end{array}$$
if $(y,t)\in\Gamma ^{p}_{\gamma}(x)$ and $$0<t<\min\left\{\log
\left(\dfrac{4k+2\delta}
{4k+\delta}\right),\log\left(\dfrac{k-1-\delta/2}{k-1-3\delta/4}
\right),\frac{1}{|x|^{2}}\wedge \frac{1}{4}\right\}=:a.$$ 
Thus
taking supremum on $(y,t)\in\Gamma ^{p}_{\gamma}(x)$,
$0<t<\alpha<a$ and then taking $\alpha\to 0^{+}$ we obtain,
$$\Omega f(x)\leq C_{d}(\epsilon+M_{\gamma_{d}}f_{2}(x))$$
for all $\epsilon>0$ and almost every  $x$ with $\mid x\mid \leq k-1.$

Given $\epsilon>0$, let us take $k$ sufficiently large such that
$$\|f_{2}\|_{1,\gamma_{d}}\leq C_{d}\epsilon^{2},$$
 then by the estimation of $\Omega$
and the weak continuity of $M_{\gamma_{d}}$ we get
$$\gamma_{d}(\{ x\in \mathbb{R}^{d}:\mid x\mid\leq k-1, \Omega f(x)>\epsilon \})\leq\epsilon$$
and that implies that
 $\Omega f(x)=0$ a.e. \ep\\

A similar proof for the Poisson-Hermite semigroup, using the non-tangential maximal function defined as
\begin{equation}
\mathcal{P}^{\ast}_{\gamma}f(x)=\sup_{(y,t)\in\Gamma_{\gamma}(x)}\mid
P_{t}f(y)\mid,
\end{equation}
and its analogous truncated version, should be possible but it has some technical difficulties that we have been unable to overcome so far.\\

\section{Non-tangential convergence of the
Ornstein-Uhlenbeck semigroup: alternative proof.}

Let us now prove a general statement for families of linear
operators that will allow us to get a simpler proof of the
non-tangential convergence, both for the Ornstein-Uhlenbeck semigroup and  for
the Poisson-Hermite semigroup. It is a generalization of Theorem
2.2 of J. Duoandikoetxea's book \cite{duo}.

\begin{thm} Let  $\{T_{t}\}_{t>0}$ be a family of  linear
operators on $L^{p}(\mathbb{R}^{d},\mu)$ and for any $x\in
\mathbb{R}^{d}$, let $\Gamma(x)$ be a subset of
$\mathbb{R}_{+}^{d+1}$ such that $x$ is in $(\Gamma(x))'$, that is to say $x$ is an 
accumulation point of $\Gamma(x)$. Let us
define
\begin{equation*}
T^{*}f(x)=\sup\{|T_{t}f(y)|:(y,t)\in\Gamma(x)\},
\end{equation*}
for $f\in  L^{p}(\mathbb{R}^{d},\mu)$ and $x\in
\mathbb{R}^{d}.$
If $T^{*}$ is weak $(p,q)$ then the set
 $$S=\left\{ f\in
L^{p}(\mathbb{R}^{d},\mu):\lim_{(y,t)\rightarrow x, (y,t)\in\Gamma(x)}T_{t}f(y)=f(x)
\, a.e. \,\right\}$$ is closed in $L^{p}(\mathbb{R}^{d},\mu)$.
\end{thm}

\dem Let us consider a equence $(f_{n})$ in $S$ such that
$f_{n}\rightarrow f$ en $L^{p}(\mathbb{R}^{d},\mu)$, then
$$|T_{t}f(y)-f(x)|-|T_{t}f_{n}(y)-f_{n}(x)|\leq
|T_{t}(f-f_{n})(y)-(f(x)-f_{n}(x))|,$$ this implies that for each
$n\in \mathbb{N}$, for almost every $x$,

\begin{eqnarray*}
\limsup_{(y,t)\rightarrow x, (y,t)\in\Gamma(x)}&&|T_{t}f(y)-f(x)| \\& \leq &
\limsup_{(y,t)\rightarrow x, (y,t)\in\Gamma(x)}|T_{t}(f-f_{n})(y)-(f(x)-f_{n}(x))| \\
& \leq & \limsup_{(y,t)\rightarrow x, (y,t)\in\Gamma(x)}|T_{t}(f-f_{n})(y)|\\
&&\quad \quad \quad \quad +
\limsup_{(y,t)\rightarrow x, (y,t)\in\Gamma(x)}|f(x)-f_{n}(x)|   \\
& \leq &
T^{*}(f-f_{n})(x)+|f(x)-f_{n}(x)|.
\end{eqnarray*}
On the other hand, if  we know that $a\leq b+c$ then $a>\lambda$
implies $b>\frac{\lambda}{2}\vee c>\frac{\lambda}{2}$.

Then, given  $\lambda >0$ and $n\in \mathbb{N}$,
 $\limsup_{(y,t)\rightarrow x, (y,t)\in\Gamma(x)}|T_{t}f(y)-f(x)| >\lambda$ implies
$$T^{*}(f-f_{n})(x)> \frac{\lambda}{2} \vee |f(x)-f_{n}(x)|  >
\frac{\lambda}{2}  \, a.e. $$
and this implies that, given  $\lambda >0$,
\begin{eqnarray*}
&&\mu\left(\left\{x:\limsup_{(y,t)\rightarrow x, (y,t)\in\Gamma(x)}|T_{t}f(y)-f(x)|
>\lambda \right\}\right) \\& &\quad \quad \quad \quad \leq 
\mu\left(\left\{x:T^{*}(f-f_{n})(x)>\frac{\lambda}{2}\right\}\right) \\
& & \quad \quad \quad \quad \quad \quad +\mu\left(\left\{x:|f(x)-f_{n}(x)|>\frac{\lambda}{2}\right\}\right)\\
\end{eqnarray*}
$$
\quad \quad \quad \quad\leq \left(\frac{2C}{\lambda}\|f-f_{n}\|_{p}\big)^{q}+\big(\frac{2}{\lambda}\|f-f_{n}\|_{p}\right)^{p},$$

\noindent for all $n\in \mathbb{N}$. Therefore,
$$\mu\left(\left\{x:\limsup_{(y,t)\rightarrow x, (y,t)\in\Gamma(x)}|T_{t}f(y)-f(x)|
>\lambda\right\}\right)=0$$
and since  this is true for all $\lambda >0$, we get that
$$\mu\left(\left\{x:\limsup_{(y,t)\rightarrow x, (y,t)\in\Gamma(x)}|T_{t}f(y)-f(x)|
>0\right\}\right)=0,$$
as 
\begin{eqnarray*}
&&\left\{x:\limsup_{(y,t)\rightarrow x, (y,t)\in\Gamma(x)}|T_{t}f(y)-f(x)|
>0\right\}\\
&&\quad \quad \quad \quad \quad \quad \quad \quad=\bigcup_{n=1}^{\infty}\left\{x:\limsup_{(y,t)\rightarrow x, (y,t)\in\Gamma(x)}|T_{t}f(y)-f(x)|
>\frac{1}{n}\right\}.
\end{eqnarray*}

Thus
$$\lim_{(y,t)\rightarrow x, (y,t)\in\Gamma(x)}T_{t}f(y)=f(x) \, a.e.$$
and then $f\in S$. Therefore $S$ is a closed set in
$L^{p}(\mathbb{R}^{d},\mu)$.\ep \\

Finally, as a consequence of this result, we get  the non-tangential convergence for
the Ornstein-Uhlenbeck semigroup $\{T_{t}\}_{t>0}$ and
the Poisson-Hermite semigroup $\{P_{t}\}_{t>0}$.

\begin{cor} The   Ornstein-Uhlenbeck semigroup $\{T_{t}\}_{t>0}$ and
the Poisson-Hermite semigroup $\{P_{t}\}_{t>0}$
verify
 $$\lim_{(y,t)\rightarrow x, (y,t) \in \Gamma^{p}_{\gamma}(x)}T_{t}f(y)=f(x) \, a.e. \, x ,$$
 $$\lim_{(y,t)\rightarrow x, (y,t) \in \Gamma_{\gamma}(x)}P_{t}f(y)=f(x)
\, a.e. \, x.$$
\end{cor}

\dem Let us discuss the proof for the the   Ornstein-Uhlenbeck semigroup $\{T_{t}\}_{t>0}$. The
proof for the  Poisson-Hermite semigroup $\{P_{t}\}_{t>0}$ is totally similar.

 It is immediate that for any given polynomial $f(x)=\sum_{k=0}^{n}J_{k}f(x)$, since 
$T_{t}f(y)=T_{t}\big(\sum_{k=0}^{n}J_{k}f(y)\big)=\sum_{k=0}^{n}e^{-tk}J_{k}f(y)$, we have
the non-tangential convergence,
 $$\lim_{(y,t)\rightarrow x, (y,t) \in \Gamma^{p}_{\gamma}(x)}T_{t}f(y)=f(x),$$
for all $x \in \mathbb{R}^{d}$. Now considering the set
 $$S=\left\{ f\in
L^{p}(\gamma_d):\lim_{(y,t)\rightarrow x, (y,t)\in\Gamma^{p}_{\gamma}(x)}T_{t}f(y)=f(x)
\, a.e. \,\right\},$$
corresponding to the Ornstein-Uhlenbeck semigroup, then the polynomials are in $S$. From the previous result, since  non-tangential maximal function for the Ornstein-Uhlenbeck semigroup $\mathcal{T}^{\ast}_{\gamma}f$ is weak $(1,1)$  with respect to the Gaussian measure, 
we get that the set 
$S$ is closed in $L^p(\gamma_d)$ and since the polynomials are dense in  $L^p(\gamma_d)$
then $S=L^{p}(\gamma_d)$. \ep\\

We want to thank the referees for their suggestions and/or corrections that improved the presentation of this paper.

\end{document}